\documentclass[10pt]{amsart}
\usepackage{latexsym, amsmath}

\newcommand{\rt}{\mathbb{R}^3}

\renewcommand{\epsilon}{\varepsilon}
\newcommand{\newsection}[1]
{\subsection{#1}\setcounter{theorem}{0} \setcounter{equation}{0}
\par\noindent}

\newtheorem{theorem}{Theorem}

\newtheorem{lemma}[theorem]{Lemma}
\newtheorem{corr}[theorem]{Corollary}

\newtheorem{proposition}[theorem]{Proposition}
\newtheorem{deff}[theorem]{Definition}

\newcommand{\bth}{\begin{theorem}}
\newcommand{\ble}{\begin{lemma}}
\newcommand{\bcor}{\begin{corr}}

\newcommand{\bdeff}{\begin{deff}}

\newcommand{\bprop}{\begin{proposition}}
\newcommand{\eth}{\end{theorem}}
\newcommand{\ele}{\end{lemma}}
\newcommand{\ecor}{\end{corr}}
\newcommand{\edeff}{\end{deff}}

\newcommand{\eprop}{\end{proposition}}

\newcommand{\cd}{\, \cdot\, }

\newcommand{\Rn}{{\mathbb R}^n}

\renewcommand{\Pi}{\varPi}

\renewcommand{\epsilon}{\varepsilon}

\newcommand{\Rplus}{{\Bbb R}_+}

\newcommand{\tidle}{\tilde}

\newcommand{\obstacle}{\mathcal{K}}
\newcommand{\varespilon}{\varepsilon}

\newcommand{\ext}{{\R^3\backslash\mathcal{K}}}

\newcommand{\R}{{\mathbb R}}

\begin{document}

\title[Dirichlet-wave equation]
{Estimates for the Dirichlet-wave equation and applications to
nonlinear wave equations}
\thanks{The author was supported in part by the NSF}

\author{Christopher D. Sogge}
\address{Department of Mathematics,  Johns Hopkins University,
Baltimore, MD 21218}

\maketitle

\newsection{Introduction}

In this  article we shall go over recent work in proving
dispersive and Strichartz estimates for the Dirichlet-wave
equation.  We shall discuss applications to existence questions
outside of obstacles and discuss open problems.

The estimates that we shall discuss involve solutions of the
Dirichlet-wave equation outside of a fixed obstacle
$\mathcal{K}\subset\mathbb{R}^n$, i.e., if $\square =
\partial_t^2-\Delta$,
\begin{equation}\label{1.1}
\begin{cases} \square u(t,x)=F(t,x), \quad t >0, \quad x\in
\Rn\backslash \obstacle
\\
u(t,x)=0, \quad t>0, \quad x\in \partial\mathcal{K}
\\
u(0,x)=f(x), \quad \partial_tu(0,x)=g(x).
\end{cases}
\end{equation}
We shall assume throughout that $\mathcal{K}$ has $C^\infty$
boundary. We also shall assume that $\mathcal{K}$ is compact, and,
by rescaling, there is no loss of generality in assuming in what
follows that
$$\mathcal{K}\subset \{x\in \mathbb{R}^n: \, |x|<1\}.$$

We shall mainly concern ourselves with the physically important
case where the spatial dimension $n$ equals $3$.  It is
considerably easier to prove estimates for the wave equation in
odd-spatial dimensions in part because of the fact that  the sharp
Huygens principle holds in this case for solutions of the
boundaryless wave equation in Minkowski space $\mathbb{R}_+\times
\Rn$.  By this we mean that if $v$ solves the Minkowski wave
equation $\square v(t,x)=0$ and if its initial data $(v(0,\cd),
\partial_tv(0,\cd))$ vanish when $|x|>R$, then $v(t,x)=0$ if $|\,
t-|x|\, |>R$.

Sharp Huygens principle of course does not hold for the obstacle
case \eqref{1.1}.  On the other hand, for a wide class of
obstacles, there is exponential decay of local energies for
compactly supported data when the spatial dimension $n$ is odd.
Specifically, in this case, if $\obstacle\subset \Rn$ is
nontrapping and if $v$ solves the homogeneous Dirichlet-wave
equation
\begin{equation}\label{1.2}
\begin{cases}
\square v(t,x)=0, \quad t>0, \quad x\in \Rn\backslash \obstacle
\\
v(t,x)=0, \quad t>0, \quad x\in \partial\obstacle,
\end{cases}
\end{equation}
then there is a constant $c>0$ so that if $R>1$ is fixed and if
\begin{equation}\label{1.3}
v(0,x)=\partial_t v(0,x)=0, \quad \{x\in \Rn\backslash \obstacle:
\, |x|>R\},\end{equation} then
\begin{equation}\label{1.4}
\Bigl(\int_{|x|<R} |v'(t,x)|^2\, dx\Bigr)^{1/2} \le Ce^{-ct}
\|v'(0,\cd)\|_2.
\end{equation}
Here, and in what follows,
$$v'=(\partial_t v,\nabla_x v)$$
denotes the space-time gradient of $v$ and in the obstacle case
the region $\{|x|<R\}$ is understood to mean $\{x\in \Rn\backslash
\obstacle: \, |x|<R\}$.

The exponential local decay of energies for nontrapping obstacles
in odd dimensions is due to Morawetz, Ralston and Strauss
\cite{MRS}, following earlier work for star-shaped obstacles of
Lax, Morawetz and Phillips \cite{LMP}.  Estimate \eqref{1.4} will
be a substitute for sharp Huygens principle that will allow us, in
certain cases, to prove global estimates, such as Strichartz
estimates, if local in time estimates hold for the obstacle case
and if the corresponding global estimates hold for Minkowski
space.

By using the local exponential decay of energy we can prove the
following sharp weighted space-time estimate for solutions of
\eqref{1.1}
\begin{multline}\label{1.5}
\Bigl(\log (2+T)\Bigr)^{-1/2}\bigl\| \, (1+|x|)^{-1/2}u'
\bigr\|_{L^2(\{(t,x)\in [0,T]\times \Rn \backslash \obstacle\})}
\\
\le C\|u'(0,\cd)\|_2 + C\int_0^T\|F(t,\cd)\|_2 \, dt,
\end{multline}
if $\obstacle$ is non-trapping and $n$ is odd.  In the region
where $|x|$ is small compared to $t$, say $|x|<t/2$, this estimate
is in some ways stronger than the usual energy estimate.  For this
reason, it plays an important role in applications to nonlinear
problems involving obstacles.  One uses \eqref{1.5} to handle
various local terms near the boundary that arise in the proofs of
the main pointwise and $L^2$ estimates.

Even though \eqref{1.4} cannot hold if there are trapped rays a
weaker form of this inequality is valid when $n$ is odd in certain
situations where there are elliptic trapped  rays.  Indeed, a
remarkable result of Ikawa \cite{Ikawa1}, \cite{Ikawa2} says that
if $v$ solves \eqref{1.2} and if \eqref{1.3} holds then
\begin{equation}\label{1.6}
\|v'(t,\cd)\|_{L^2(|x|<R)}\le C e^{-ct}\sum_{|\alpha|\le
1}\|\partial_x^\alpha v'(0,\cd)\|_{L^2(|x|<R)},
\end{equation}
for some constant $c>0$ if $\obstacle$ is a finite union of convex
obstacles. In the case of three or more obstacles Ikawa's result
requires a technical assumption that the obstacles are
sufficiently separated, but it is thought that \eqref{1.6} should
hold in the case where there are no hyperbolic trapped rays. Also,
just by interpolating with the standard energy estimate, one
concludes that the variant of \eqref{1.6} holds if one replaces
the $L^2$ norm of $v'(0,\cd)$ by an $H^{\varepsilon}$ norm with
$\varepsilon>0$ and the constant $c>0$ in the exponential
depending on $\varespilon$.  This fact would allow one to prove
global Strichartz estimates with arbitrary small loss of
derivatives if the local in time estimates were known (cf.
\cite{burq3}).  For other local decay bounds see Burq \cite{burq}.

In the rest of the paper we shall indicate how one can use the
exponential local decay of energy to prove global estimates for
solutions of \eqref{1.1} that have applications to nonlinear
Dirichlet-wave equations.  In the next section we shall go over
the simplest situation of proving global Strichartz estimates in
$\rt \backslash \obstacle$ when $\obstacle$ is  convex with smooth
boundary.  This argument  will serve as a template for the more
involved ones that are used to prove almost global and global
existence for certain quasilinear wave equations.  The most basic
of these, which will be discussed in \S 3, will be to show that
one can prove fixed-time $L^2$ estimates and weighted space-time
$L^2$ estimates for $\Omega_{ij}u'$ if
\begin{equation}\label{1.7}
\Omega_{ij}=x_i\partial_j-x_j\partial_i, \quad 1\le i<j\le 3,
\end{equation}
 are angular-momentum operators for $\rt$.  As we shall see, by using these
estimates one can prove almost global existence for semilinear
wave equations in $\rt \backslash \obstacle$ if $\obstacle$ is
nontrapping.  In the next section we shall see how one can prove a
pointwise dispersive estimate for solutions of \eqref{1.1} if
$\obstacle$ is nontrapping or if it satisfies Ikawa's conditions.
We shall also present related $L^2$ estimates that can be used to
prove almost global existence results for quasilinear
Dirichlet-wave equations and global existence for ones satisfying
an appropriate null condition.

The results described in this paper were presented in a series of
lectures given by the author in Japan in July of 2002.  The author
is grateful for the hospitality shown to him, especially that of
H. Kozono and M. Yamazaki.

\newsection{Strichartz estimates outside convex obstacles}

In this section we shall show how local Strichartz estimates for
obstacles, global ones for Minkowski space and the energy decay
estimates \eqref{1.4} can be used to prove global Strichartz
estimates for obstacles.  This was first done in the case of odd
dimensions by Smith and the author \cite{SS}, and later for even
dimension by Burq \cite{burq2} and Metcalfe \cite{Met}.

For simplicity, we shall only consider the special case where the
spatial dimension $n$ is equal to three.  We shall also only treat
the most basic Strichartz estimate in this case.  The global
Minkowski version, which will be used in the proof of the version
for obstacles, says that
\begin{equation}\label{2.1}
\|v\|_{L^4({\mathbb R}_+\times \rt)}\le C\Bigl(\|v(0,\cd)\|_{\dot
H^{1/2}(\rt)}+\|\partial_t v(0,\cd)\|_{\dot H^{-1/2}(\rt)}+
\|\square v\|_{L^{4/3}({\mathbb R}_+\times \rt)}\Bigr).
\end{equation}
Here $\dot H^\gamma(\rt)$ denote the homogeneous Sobolev spaces on
$\rt$.

In addition to this, if $\obstacle \subset \rt$ is our compact
obstacle, we shall need to assume that we have the local in time
Strichartz estimates
\begin{multline}\label{2.2}
\|u\|_{L^4([0,1]\times \rt\backslash \obstacle)} + \sup_{0\le t\le
1}\Bigl(\|u(t,\cd)\|_{H^{1/2}_D(\rt\backslash
\obstacle)}+\|u(t,\cd)\|_{H^{-1/2}_D(\rt\backslash
\obstacle)}\Bigr)
\\
\le C\Bigl(\|u(0,\cd)\|_{H^{1/2}_D(\rt\backslash
\obstacle)}+\|\partial_t u(0,\cd)\|_{H^{-1/2}_D(\rt\backslash
\obstacle)} + \|F\|_{L^{4/3}([0,1]\times \rt\backslash
\obstacle)}\Bigr),
\end{multline}
assuming that the initial data is supported in the set $\{x\in \rt
\backslash \obstacle: \, |x|<4\}$. Here, $H^\gamma_D(\rt\backslash
\obstacle)$ are the usual Dirichlet-Sobolev spaces.

For the homogeneous case where the forcing term $F\equiv 0$ it was
shown in \cite{SS2} that \eqref{2.2} holds when $\obstacle \subset
\rt$ is convex.  An interesting problem would be to show that this
estimate holds for a larger class of obstacles.  In \cite{SS2}
more general Strichartz estimates for convex obstacles in all
dimensions were also proved.  In \cite{SS} estimates for the
inhomogeneous wave equation were also obtained by using a lemma of
Christ and Kiselev \cite{CK}.

In addition to \eqref{2.1} and \eqref{2.2}, we shall need a
Sobolev space variant of \eqref{1.4}.  We suppose that $R>1$ is
given and that $\beta(x)$ is smooth and supported in $|x|\le R$.
Then there is a $c>0$ so that
\begin{equation}\label{2.3}
\|\beta u(t,\cd)\|_{H^{1/2}_D}+\|\beta \partial_t
u(t,\cd)\|_{H^{-1/2}_D} \le
Ce^{-ct}\Bigl(\|u(0,\cd)\|_{H^{1/2}_D}+\|\partial_t
u(0,\cd)\|_{H^{-1/2}_D}\Bigr),
\end{equation}
if $u$ solves \eqref{1.1} with vanishing forcing term $F$ and has
initial data satisfying $u(0,x)=\partial_tu(0,x)=0$, $|x|>R$. This
estimate just follows from \eqref{1.4} and a simple interpolation
argument.

We claim that by using these three inequalities, we can prove the
following result from \cite{SS}.

\begin{theorem}\label{theorem2.1}  Let $u$ solve \eqref{1.1} when $\obstacle\subset
\rt$ is a convex obstacle with smooth boundary.  Then
\begin{equation}\label{2.4}
\|u\|_{L^4({\mathbb R}_+\times \rt\backslash \obstacle)}\le C
\Bigl(
\|f\|_{H^{1/2}_D}+\|g\|_{H^{-1/2}_D}+\|F\|_{L^{4/3}({\mathbb
R}_+\times \rt\backslash\obstacle)}\Bigr).
\end{equation}
\end{theorem}

Recall that we are assuming, as we may, that $\obstacle \subset
\{x\in \rt: \, |x|<1\}$.  The first step in the proof of this
result will be to establish the following

\begin{lemma}\label{lemma2.2}  Let $u$ solve the Cauchy problem
\eqref{1.1} with forcing term $F$ replaced by $F+G$.  Suppose that
the initial data is supported in $\{\, |x|\le 2\}$ and that $F$,
$G$ are supported in $\{ 0\le t\le 1\} \times \{\, |x|\le 2\}$.
Then if $\rho<c$, where $c$ is the constant in \eqref{2.3},
\begin{multline}\label{2.5}
\Bigl\| e^{\rho(t-|x|)}u\|_{L^4(\mathbb{R}_+\times \rt\backslash
\obstacle)}
\\
\le
C\Bigl(\|f\|_{H^{1/2}_D}+\|g\|_{H^{-1/2}_D}+\|F\|_{L^{4/3}(\mathbb{R}_+\times
\rt\backslash \obstacle)}+\int
\|G(t,\cd)\|_{H^{-1/2}_D(\rt\backslash \obstacle)}\, dt\Bigr).
\end{multline}
\end{lemma}

\noindent {\bf Proof of Lemma \ref{lemma2.2}:}  By \eqref{2.2} and
Duhamel's principle, the inequality holds for the $L^4(dtdx)$ norm
of $u$ over $[0,1]\times \rt\backslash \obstacle$.  Also, by
\eqref{2.2},
\begin{multline}
\|u(1,\cd)\|_{H^{1/2}_D(\rt\backslash
\obstacle)}+\|\partial_tu(1,\cd)\|_{H^{-1/2}_D(\rt\backslash
\obstacle)}
\\
\le
C\Bigl(\|f\|_{H^{1/2}_D}+\|g\|_{H^{-1/2}_D}+\|F\|_{L^{4/3}({\mathbb
R}_+\times \rt\backslash \obstacle)}+\int
\|G(t,\cd)\|_{H^{-1/2}_D(\rt\backslash \obstacle)}\, dt\Bigr).
\end{multline}
By considering $t\ge 1$, we may take $F=G=0$, with $(f,g)$ now
supported in $\{\, |x|\le 3\}$.

We next decompose $u=\beta u+(1-\beta)u$, where $\beta(x)=1$ for
$|x|\le 1$ and $\beta(x)=0$ for $|x|\ge2$.  Let us first consider
$\beta u$.  We write
$$(\partial_t^2-\Delta)(\beta u) = -2\nabla_x \beta\cdot \nabla_x
u -(\Delta \beta) u = \tilde G(t,x),$$ and note that $\tidle
G(t,x)=0$ if $|x|\ge 2$.  By \eqref{2.3} we have
\begin{multline}\label{2.7}
\|\tilde G(t,\cd)\|_{H^{-1/2}_D}+\|\beta
u(t,\cd)\|_{H^{1/2}_D}+\|\partial_t (\beta
u)(t,\cd)\|_{H^{-1/2}_D}
\\
\le C e^{-ct}\Bigl(\|f\|_{H^{1/2}_D}+\|g\|_{H^{-1/2}_D}\Bigr).
\end{multline}
By \eqref{2.2} and Duhamel's principle, it follows that
$$\|\beta u\|_{L^4([j,j+1]\times \rt\backslash \obstacle)}\le
Ce^{-cj}\Bigl(\|f\|_{H^{1/2}_D}+\|g\|_{H^{-1/2}_D}\Bigr),$$ which
implies that $\beta u$ satisfies the bounds in \eqref{2.5}.

Now let us show that the same is true for $(1-\beta)u$.  On the
support of $(1-\beta)u$, we have
$$(\partial_t^2-\Delta)u=-\tilde G,$$
and by Duhamel's principle we have
$$u(t,x)=u_0(t,x)+\int_0^t u_s(t,x)\, ds,$$
where $u_0$ is the solution of the Minkowski wave equation on
${\mathbb R}_+\times \rt$ with initial data $\bigl(\, (1-\beta)f,
(1-\beta)g\, \bigr)$, and where $u_s(t,x)$ is the solution of the
Minkowski space wave equation on the set $t>s$ with Cauchy data
$(0,\tilde G(s,\cd))$ on the hyperplane $t=s$.  (Recall that
$\tilde G$ and $(1-\beta)$ vanish near $\partial \obstacle$.)
Since the initial data of $u_0$ is supported in $\{x\in \rt: \,
|x|\le 2\}$, by the sharp Huygens principle, $u_0$ must satisfy
the bounds in \eqref{2.5}.  Additionally, on the support of
$u_s(t,x)$ have $t\ge s$ and $t-|x|\in [s-3,s+3]$, so that by
\eqref{2.1} and \eqref{2.7} we have
$$\bigl\| e^{\rho( t-|x|)}u_s\|_{L^4(dtdx)}\le
Ce^{(\rho-c)s}\Bigl(\|f\|_{H^{1/2}_D}+\|g\|_{H^{-1/2}_D}\Bigr),$$
which leads to the desired estimate for the remaining part of $u$.
\qed

We also require a simple consequence of Plancherel's theorem:

\begin{lemma}\label{lemma2.3}  Let $\beta(x)$ be smooth and
supported in $\{x\in \rt: \, |x|\le 2\}$.  Then
$$\int_{-\infty}^{+\infty}\|\beta(\cd)
\Bigl(e^{it|D|}f\Bigr)(t,\cd)\|^2_{H^{1/2}(\rt)}\, dt \le
C\|f\|_{\dot H^{1/2}(\rt)}^2,$$ if $|D|=\sqrt{-\Delta}$.
\end{lemma}

\noindent{\bf Proof:}  By Plancherel's theorem over $t,x$, the
left side can be written as
$$\int_0^\infty\int \Bigl| \int \hat \beta(\xi-\eta)\hat
f(\eta)\delta(\tau-|\eta|)\, d\eta\Bigr|^2 \, (1+|\xi|^2)^{1/2} \,
d\xi d\tau.$$ If we apply the Schwarz inequality in $\eta$ we
conclude that this is dominated by
\begin{multline*}
\int_0^\infty \int \Bigl( \int|\hat\beta(\xi-\eta)|
\delta(\tau-|\eta|) \, d\eta \Bigr) \Bigl(\int |\hat
\beta(\xi-\eta)| \, |\hat f(\eta)|^2 \delta(\tau-|\eta|)\,
d\eta\Bigr)
\\
\times (1+|\xi|^2)^{1/2}\, d\xi d\tau.
\end{multline*}
This in turn is dominated by
$$\int |\hat f(\eta)|^2 \min\Bigl( |\eta|^{2},
(1+|\eta|^2)^{1/2}\Bigr) \, d\eta \le C\|f\|^2_{\dot
H^{1/2}(\rt)},$$ since
$$\sup_\xi (1+|\xi|^2)^{1/2}\Bigl( \int |\hat \beta(\xi-\eta)|\,
\delta(\tau-|\eta|) \, d\eta \Bigr) \le C\min \bigl(\tau^2,
(1+\tau^2)^{1/2}\bigr),$$ which completes the proof. \qed

\begin{corr}\label{corollary2.4}  Let $\beta$ be as above,
and let $u$ solve the ${\mathbb R}_+\times \rt$ Minkowski wave
equation $\square u=F$ with initial data $(f,g)$. Then
$$\sum_{|\alpha|\le 1}\int_0^\infty \|\beta
\partial^\alpha_{t,x}u(t,\cd)\|_{H^{-1/2}(\rt)}^2 \, dt
\le C\Bigl(\|f\|_{\dot H^{1/2}}+\|g\|_{\dot
H^{-1/2}}+\|F\|_{L^{4/3}({\mathbb R}_+\times \rt)} \Bigr)^2.$$
\end{corr}

\noindent{\bf Proof:}  If $F=0$ then this is a direct consequence
of the preceding lemma.  If $f=g=0$ then the Minkowski Strichartz
estimate \eqref{2.1}, duality, and Huygens principle imply that
for $t>0$
$$\sum_{|\alpha|\le 1}\|\beta
\partial^\alpha_{t,x}u(t,\cd)\|^2_{H^{-1/2}(\rt)}\le
C\|F\|^2_{L^{4/3}(\Gamma_t)},$$ where
$$\Gamma_t=\{(s,x): \, s\ge 0, \, s+|x|\in [t-2,t+2]\}.$$
Since $4/3\le 2$,
$$\int_0^\infty\|F\|^2_{L^{4/3}(\Gamma_t)}\, dt \le
4\|F\|^2_{L^{4/3}({\mathbb R}_+\times \rt)},$$ which finishes the
proof. \qed

\noindent{\bf Proof of Theorem 2.1:}  By Lemma 2.2, we may without
loss of generality assume that $f$ and $g$ vanish for $|x|\le 2$.
If $\beta$ is as above write
$$u=u_0-v=(1-\beta)u_0+\beta u_0-v,$$
where $u_0$ solves the Cauchy problem for the Minkowski wave
equation, with data $f,g,F$, where we set $F=0$ in ${\mathbb
R}_+\times \obstacle$.  By \eqref{2.1}, $u_0$ satisfies the
desired bounds, and so we just need to estimate $\beta u_0-v$.  We
write
$$(\partial_t^2-\Delta)(\beta u_0-v)=\beta F+G,$$
where $G=-2\nabla_x \beta \cdot \nabla_x u_0-(\Delta \beta)u_0$
vanishes for $|x|\ge 2$, and satisfies
\begin{equation}\label{2.8}\int_0^\infty \|G(t,\cd)\|^2_{H_D^{-1/2}}\, dt \le
C\Bigl(\|f\|_{\dot H^{1/2}}+\|g\|_{\dot H^{-1/2}}
+\|F\|_{L^{4/3}}\Bigr)^2,
\end{equation}
by Corollary \ref{corollary2.4}.  Note that the initial data of
$\beta u_0-v$ vanishes.  Let $F_j$, $G_j$ denote the restricitions
of $F,G$ to the set where $t\in [j,j+1]$, and write for $t>0$
$$\beta u_0-v =\sum_{j=0}^\infty u_j(t,x),$$
where $u_j$ is the forward solution of $(\partial_t^2-\Delta)u_j =
\beta F_j+G_j$.

By Lemma \ref{lemma2.2}, the following holds
$$\Bigl\| \, e^{\rho(t-j-|x|)} u_j \, \Bigr\|_{L^4}\le C\Bigl(\|\beta
F_j\|_{L^{4/3}}+\int_j^{j+1}\|G(t,\cd)\|_{H^{-1/2}}\, dt\Bigr).$$
Furthermore, $u_j(t,x)$ is supported on the set where $t-j-|x|\ge
-2$.  Consequently, we have
\begin{align*}
\|\beta u_0-v\|_{L^4(dtdx)}^2 &\le C\sum_{j=0}^\infty
\|e^{\rho(t-j-|x|)}u_j\|^2_{L^4(dtdx)}
\\
&\le C\sum_{j=0}^\infty
\|F_j\|^2_{L^{4/3}}+C\sum_{j=0}^\infty\Bigl(\int_j^{j+1}\|G(t,\cd)\|_{H^{-1/2}}\,
dt\Bigr)^2
\\
&\le C\|F\|^2_{L^{4/3}} +C\int_0^\infty\|G(t,\cd)\|^2_{H^{-1/2}}\,
dt.
\end{align*}
If we use \eqref{2.8}, we conclude that $\beta u_0-v$ also
satisfies the desired bounds, which completes the proof. \qed

\noindent {\bf Remark:}  It would be very interesting to see
whether the Strichartz estimates of Georgiev, Lindblad and the
author \cite{gls} or Tataru \cite{T} are valid for ${\mathbb
R}_+\times \rt \backslash \obstacle$ when, as above $\obstacle$ is
convex.

\newsection{Weighted space-time $L^2$ estimates}

In \cite{KSS2}, the following weighted space-time estimate for
Minkowski space was proved
\begin{multline}\label{3.1}
\bigl(\log (2+T)\bigr)^{-1/2} \|(1+|x|)^{-1/2}v'\|_{L^2(\{(t,x):
0\le t\le T, \, x\in \rt\})}
\\
\le C\|v'(0,\cd)\|_{L^2(\rt)}+C\int_0^T \|\square
v(t,\cd)\|_{L^2(\rt)}\, dt.
\end{multline}
By using this estimate and the exponential local decay of energy,
one can adapt the arguments of the previous section to prove the
following analogous estimates for solutions of the Dirichlet-wave
equation \eqref{1.1} if $\obstacle \subset \rt$ is non-trapping
\begin{multline}\label{3.2}
\Bigl(\log (2+T)\Bigr)^{-1/2}\bigl\| \, (1+|x|)^{-1/2}u'
\bigr\|_{L^2(\{(t,x)\in [0,T]\times \rt \backslash \obstacle\})}
\\
\le C\|u'(0,\cd)\|_{L^2(\rt\backslash\obstacle)} +
C\int_0^T\|F(t,\cd)\|_{L^2(\rt\backslash \obstacle)} \, dt.
\end{multline}
Indeed, the proof of Lemma \ref{lemma2.3} shows that this estimate
is valid (without the $\log$ weight) when one replaces the $L^2$
norm in the left side of \eqref{3.2} by one over $\{(t,x): 0\le
t\le T, x\in \rt \backslash \obstacle, |x|<2\}$, with a constant
that is independent of $T$. Using this and the Minkowski space
estimates \eqref{3.1}, one sees that the analog of \eqref{3.2}
also holds when the norm is taken over the region where $|x|>2$.

To handle applications to nonlinear wave equations, one requires a
slight generalization of this estimate, which involves the
operators
\begin{equation}\label{3.3}
Z=\{\partial_t, \partial_i, \Omega_{jk}, 1\le i\le 3, 1\le j<k\le
3\}.
\end{equation}

\begin{theorem}\label{theorem3.1}  If $u$ is as in \eqref{1.1} has
vanishing Cauchy data, then for any $N=0,1,2,\dots$
\begin{multline}\label{3.4}
\sum_{|\alpha|\le N}\Bigl( \|Z^\alpha u'(t,\, \cdot\,
)\|_{L^2({\mathbb R}^3\backslash {\mathcal K})}
+(\ln(2+t))^{-1/2}\|(1+|x|)^{-1/2}Z^\alpha u'\|_{L^2(\{(s,x)\in
[0,t]\times {\mathbb R}^3\backslash {\mathcal K}\})}\Bigr)
\\
\le C\sum_{|\alpha|\le N}\int_0^t \|Z^\alpha F(s,\, \cdot\,
)\|_{L^2({\mathbb R}^3\backslash {\mathcal K})} \, ds +
 C\sup_{0\le s\le t}\sum_{|\alpha|\le N-1}\|Z^\alpha F(s,\,
\cdot\, )\|_{L^2({\mathbb R}^3\backslash {\mathcal K})}
\\
+ C\sum_{|\alpha|\le N-1}\|Z^\alpha F\|_{L^2(\{(s,x)\in
[0,t]\times {\mathbb R}^3\backslash {\mathcal K}\})}.
\end{multline}
\end{theorem}

Let us first see that the estimate holds when the norm in the left
is taken over $|x|<2$.  Clearly the first term in the left is
under control since
\begin{equation*}\sum_{|\alpha|\le
N}\|Z^\alpha u'(t,\, \cdot\, )\|_{L^2(\{x\in {\mathbb
R}^3\backslash {\mathcal K}): \, |x|<2\}}
\le C_N \sum_{|\alpha|\le N}\|\partial_{t,x}^\alpha u'(t,\,
\cdot\, )\|_{L^2({\mathbb R}^3\backslash {\mathcal
K})},\end{equation*} and standard arguments imply that the right
hand side here is dominated by
\begin{multline}\label{3.5}
\sum_{|\alpha|\le N}\|\partial_{t,x}^\alpha u'(t,\, \cdot\,
)\|_{L^2({\mathbb R}^3\backslash {\mathcal K})}
\\
\le C\sum_{|\alpha|\le N}\int_0^t \|\partial_{t,x}^\alpha F(s,\,
\cdot\, )\|_{L^2({\mathbb R}^3\backslash {\mathcal K})} \, ds +C
\sum_{|\alpha|\le N-1}\|\partial_{t,x}^\alpha F(t,\, \cdot\,
)\|_{L^2({\mathbb R}^3\backslash {\mathcal K})}.
\end{multline}
Indeed, if $N=0$, \eqref{3.5} is just the standard energy
identity. To prove that \eqref{3.5} holds for $N$, assuming that
it is valid when $N$ is replaced by $N-1$, one notes that since
$\partial_t w$ vanishes on the boundary one has
\begin{multline*}
\sum_{|\alpha|\le N-1}\|\partial_{t,x}^\alpha \partial_tu'(t,\,
\cdot\, )\|_{L^2({\mathbb R}^3\backslash {\mathcal K})}
\\
\le C \sum_{|\alpha|\le N-1}\int_0^s\|\partial_{t,x}^\alpha
\partial_s F(s,\, \cdot\, )\|_{L^2({\mathbb R}^3\backslash {\mathcal
K})} +  C \sum_{|\alpha|\le N-2}\|\partial_{t,x}^\alpha \partial_t
F(t,\, \cdot\, )\|_{L^2({\mathbb R}^3\backslash {\mathcal K})}
.\end{multline*} Since $\partial^2_t w =\Delta w + F$, we get from
this that
\begin{multline*}\sum_{|\alpha|\le
N-1}\|\partial_{t,x}^\alpha \Delta u(t,\, \cdot\,
)\|_{L^2({\mathbb R}^3\backslash {\mathcal K})}
\\
 \le C
\sum_{|\alpha|\le N}\int_0^s\|\partial_{t,x}^\alpha F(s,\, \cdot\,
)\|_{L^2({\mathbb R}^3\backslash {\mathcal K})}\, ds +
\sum_{|\alpha|\le N-1}\|\partial_{t,x}^\alpha F(t,\, \cdot\,
)\|_{L^2({\mathbb R}^3\backslash {\mathcal K})}.
\end{multline*}
By elliptic regularity, $\sum_{|\alpha|\le
N}\|\partial^\alpha_{t,x}u'(t,\, \cdot\, ) \|_{L^2({\mathbb
R}^3\backslash {\mathcal K})}$ is dominated by the left side of
the last equation, which finishes the proof of \eqref{3.4}, since
$$\sum_{|\alpha|\le N}\|\partial_{t,x}^\alpha u'(t,\, \cdot\,
)\|_{L^2({\mathbb R}^3\backslash {\mathcal K})} \le
\sum_{|\alpha|\le N}\|\partial_{x}^\alpha u'(t,\, \cdot\,
)\|_{L^2({\mathbb R}^3\backslash {\mathcal K})} +
\sum_{|\alpha|\le N-1}\|\partial_{t,x}^\alpha
\partial_tu'(t,\, \cdot\, )\|_{L^2({\mathbb R}^3\backslash {\mathcal
K})}.$$

To handle the second term on the left side  of \eqref{3.3}, again
when the left hand norm is taken over $|x| < 2$, we shall need the
following

\begin{lemma}\label{lemma3.2}  If $u$ is as in \eqref{1.1} then
for any $N=0,1,2,\dots$
\begin{multline}\label{3.6}
\sum_{|\alpha|\le N}\|\partial^\alpha_{t,x}u'\|_{L^2(\{(s,x)\in
[0,t]\times {\mathbb R}^3\backslash {\mathcal K}: \, |x|<2\})}
\\
\le C\sum_{|\alpha|\le N}\int_0^t \|\partial_{t,x}^\alpha F(s,\,
\cdot\, )\|_{L^2({\mathbb R}^3\backslash {\mathcal K})} \, ds +
C\sum_{|\alpha|\le N-1}\|\partial_{t,x}^\alpha F\|_{L^2(\{(s,x)\in
[0,t]\times {\mathbb R}^3\backslash {\mathcal K}\})}.
\end{multline}
\end{lemma}

Clearly \eqref{3.6} implies that
\begin{multline*}\sum_{|\alpha|\le N}\|Z^\alpha
u'\|_{L^2(\{(s,x)\in [0,t]\times {\mathbb R}^3\backslash {\mathcal
K}: \, |x|<2\})}
\\
\le C\sum_{|\alpha|\le N}\int_0^t \|Z^\alpha F(s,\, \cdot\,
)\|_{L^2({\mathbb R}^3\backslash {\mathcal K})} \, ds +
C\sum_{|\alpha|\le N-1}\|Z^\alpha F\|_{L^2(\{(s,x)\in [0,t]\times
{\mathbb R}^3\backslash {\mathcal K}\})},
\end{multline*}
finishing the proof that the analog of \eqref{3.4} holds where the
norms in the left are taken over $|x|<2$.

\noindent{\bf Proof of Lemma \ref{lemma3.2}:}  By the proof of
\eqref{3.5}, \eqref{3.6} follows from the special case where
$N=0$:
\begin{equation}\label{3.7}
\|u'\|_{L^2(\{(s,x)\in [0,t]\times {\mathbb R}^3\backslash
{\mathcal K}: \, |x|<2\})}
\le C\int_0^t \| F(s,\, \cdot\, )\|_{L^2({\mathbb R}^3\backslash
{\mathcal K})} \, ds,
\end{equation}
which, as we noted before, follows from the proof of Lemma
\ref{lemma2.3}. \qed

\medskip

\noindent{\bf End of proof of Theorem \ref{theorem3.1}:}  We need
to see that
\begin{multline}\label{3.8}
\sum_{|\alpha|\le N}\Bigl( \|Z^\alpha u'(t,\, \cdot\,
)\|_{L^2(|x|>2)} +(\ln(2+t))^{-1/2}\|(1+|x|)^{-1/2}Z^\alpha
u'\|_{L^2(\{[0,t]\times \{x: \, |x|>2\}\})}\Bigr)
\\
\le C\sum_{|\alpha|\le N}\int_0^t \|Z^\alpha F(s,\, \cdot\,
)\|_{L^2({\mathbb R}^3\backslash {\mathcal K})} \, ds +
 C\sup_{0\le s\le t}\sum_{|\alpha|\le N-1}\|Z^\alpha F(s,\,
\cdot\, )\|_{L^2({\mathbb R}^3\backslash {\mathcal K})}
\\
+ C\sum_{|\alpha|\le N-1}\|Z^\alpha F\|_{L^2(\{(s,x)\in
[0,t]\times {\mathbb R}^3\backslash {\mathcal K}\})}.
\end{multline}

For this we fix $\beta\in C^\infty({\mathbb R}^3)$ satisfying
$\beta(x)=1$, $|x|\ge 2$ and $\beta(x)=0$, $|x|\le 3/2$.  Then
since, by the assumption that the obstacle is contained in the set
$|x|<1$, it follows that $v=\beta u$ solves the boundaryless wave
equation
$$\square v = \beta F -2\nabla_x \beta\cdot \nabla_x u - (\Delta \beta)u$$
with zero initial data, and satisfies $u(t,x)=v(t,x)$, $|x|\ge 2$.
If we split $v=v_1+v_2$, where $\square v_1=\beta F$, and $\square
v_2 =-2\nabla_x \beta\cdot \nabla_x u - (\Delta \beta)u$, it then
suffices to prove that
\begin{multline}\label{3.9}
\sum_{|\alpha|\le N}\Bigl( \|Z^\alpha v_2'(t,\, \cdot\,
)\|_{L^2(|x|>2)} +(\ln(2+t))^{-1/2}\|(1+|x|)^{-1/2}Z^\alpha
v_2'\|_{L^2(\{[0,t]\times \{x: \, |x|>2\}\})}\Bigr)
\\
\le C\sum_{|\alpha|\le N}\int_0^t \|Z^\alpha F(s,\, \cdot\,
)\|_{L^2({\mathbb R}^3\backslash {\mathcal K})} \, ds +
 C\sup_{0\le s\le t}\sum_{|\alpha|\le N-1}\|Z^\alpha F(s,\,
\cdot\, )\|_{L^2({\mathbb R}^3\backslash {\mathcal K})}
\\
+ C\sum_{|\alpha|\le N-1}\|Z^\alpha F\|_{L^2(\{(s,x)\in
[0,t]\times {\mathbb R}^3\backslash {\mathcal K}\})}.
\end{multline}
This is because by \eqref{3.1} we have
\begin{multline*}\sum_{|\alpha|\le N}\Bigl(
\|Z^\alpha v_1'(t,\, \cdot\, )\|_{L^2(|x|>2)}
+(\ln(2+t))^{-1/2}\|(1+|x|)^{-1/2}Z^\alpha
v_1'\|_{L^2(\{[0,t]\times \{x: \, |x|>2\}\})}\Bigr)
\\
\le C\sum_{|\alpha|\le N}\int_0^t \|Z^\alpha F(s,\, \cdot\,
)\|_{2} \, ds,
\end{multline*}
due to the fact that
$$\sum_{|\alpha|\le N}\int_0^t \|Z^\alpha (\beta F)(s,\, \cdot\,
)\|_{2} \, ds \le C \sum_{|\alpha|\le N}\int_0^t \|Z^\alpha F(s,\,
\cdot\, )\|_{2} \, ds.$$

To prove \eqref{3.9} we note that $G=-2\nabla_x \beta\cdot
\nabla_x u -(\Delta \beta)u=\square v_2$, vanishes unless
$1<|x|<2$.  To use this, fix $\chi\in C^\infty_0({\mathbb R})$
satisfying $\chi(s)=0$, $|s|>2$, and $\sum_j\chi(s-j)=1$. We then
split  $G=\sum_j G_j$, where $G_j(s,x)=\chi(s-j)G(s,x)$, and let
$v_{2,j}$ be the solution of the corresponding inhomogeneous wave
equation $\square v_{2,j}=G_j$ with zero initial data in Minkowski
space. By sharp Huygen's principle we have that $|Z^\alpha
v_2(t,x)|^2\le C\sum_j |Z^\alpha v_{2,j}(t,x)|^2$ for some uniform
constant $C$. Therefore, by \eqref{3.1} we have that the square of
the left side of \eqref{3.9} is dominated by
\begin{align*}
\sum_{|\alpha| \leq N} &\sum_j \Bigl(\int_0^t\|Z^{\alpha}G_j(s,\,
\cdot\, )\|_2 ds \Bigr)^2
\\
&\le C \sum_{|\alpha| \leq N} \|Z^{\alpha}G\|_{L^2(\{(s,x): \,
0\le s\le t, \, 1<|x|<2\})}^2
\\
&\le C\sum_{|\alpha|\le N}\|Z^\alpha u'\|_{L^2(\{(s,x): \, 0\le
s\le t, \, 1<|x|<2\})}^2 + C\sum_{|\alpha|\le N}\|Z^\alpha
u\|_{L^2(\{(s,x): \, 0\le s\le t, \, 1<|x|<2\})}^2
\\
&\le C \sum_{|\alpha|\le N}\|Z^\alpha u'\|_{L^2(\{(s,x)\in
[0,t]\times {\mathbb R}^3\backslash {\mathcal K}: |x|<2\})}^2
\\
&\le C \sum_{|\alpha|\le N}\|\partial_{t,x}^\alpha
u'\|_{L^2(\{(s,x)\in [0,t]\times {\mathbb R}^3\backslash {\mathcal
K}: |x|<2\})}^2.
\end{align*}
Consequently, \eqref{3.9} follows from \eqref{3.6}, which finishes
the proof. \qed

To handle almost global existence, in addition to \eqref{3.4}, we
need the following consequence of the Sobolev estimates for $S^2
\times [0,\infty)$
\begin{equation}\label{KSob}
\|h\|_{L^2(\{x\in \rt \backslash \obstacle: \, |x|\in [R-1,R]\}}
\le \frac{C}{R}\sum_{|\alpha|\le 2}\|Z^\alpha h\|_{L^2(\{x\in \rt
\backslash \obstacle: \, |x|\in [R-2,R+1]\})}, \quad R\ge 1.
\end{equation}

Let us conclude this section by showing how \eqref{3.4} and
\eqref{KSob} can be used to prove almost global existence of
semilinear wave equations outside of non-trapping obstacles.  We
shall consider semilinear systems of the form
\begin{equation}\label{0.2}
\begin{cases}
\square u = Q(u'), \quad (t,x)\in {\mathbb R}_+\times {\mathbb
R}^3\backslash {\mathcal K}
\\
u(t,\cdot)|_{\partial\mathcal K}=0
\\
u(0,\cdot)=f, \, \, \partial_tu(0,\cdot)=g.
\end{cases}
\end{equation}
Here
$$\square = \partial_t^2 - \Delta$$
is the D'Alembertian, with $\Delta =
\partial_1^2+\partial_2^2+\partial_3^2$ being the standard
Laplacian.  Also, $Q$ is a constant coefficient quadratic form in
$u'=(\partial_tu, \nabla_x u)$.

In the non-obstacle case we shall obtain almost global existence
for equations of the form
\begin{equation}\label{0.10}
\begin{cases}
\square u = Q(u'), (t,x)\in {\mathbb R}_+\times {\mathbb R}^3
\\
u(0,\cdot)=f, \, \, \partial_tu(0,\cdot)=g.
\end{cases}
\end{equation}

In order to solve \eqref{0.2} we must also assume that the data
satisfies the relevant compatibility conditions.  Since these are
well known (see e.g., \cite{KSS}), we shall describe them briefly.
To do so we first let $J_ku =\{\partial^\alpha_xu: \, 0\le
|\alpha|\le k\}$ denote the collection of all spatial derivatives
of $u$ of order up to $k$.  Then if $m$ is fixed and if $u$ is a
formal $H^m$ solution of \eqref{0.2} we can write
$\partial_t^ku(0,\cdot)=\psi_k(J_kf,J_{k-1}g)$, $0\le k\le m$, for
certain compatibility functions $\psi_k$ which depend on the
nonlinear term $Q$ as well as $J_kf$ and $J_{k-1}g$.  Having done
this, the compatibility condition for \eqref{0.2} with $(f,g)\in
H^m\times H^{m-1}$ is just the requirement that the $\psi_k$,
vanish on $\partial{\mathcal K}$ when $0\le k\le m-1$.
Additionally, we shall say that $(f,g)\in C^\infty$ satisfy the
compatibility conditions to infinite order if this condition holds
for all $m$.

If $\{\Omega\}$ denotes the collection of vector fields
$x_i\partial_j-x_j\partial_i$, $1\le i<j\le 3$, then we can now
state our existence theorem.

\begin{theorem}\label{theorem0.1}
Let ${\mathcal K}$ be a smooth compact nontrapping obstacle  and
assume that $Q(u')$ is above. Assume further that $(f,g)\in
C^\infty({\mathbb R}^3\backslash {\mathcal K})$ satisfies the
compatibility conditions to infinite order.  Then there  are
constants $c, \epsilon_0 >0$ so that if $\varepsilon \leq
\epsilon_0$ and
\begin{equation}\label{dataob}
\sum_{|\alpha|+j\le 10}\|\partial_x^j\Omega^\alpha
f\|_{L^2({\mathbb R}^3\backslash {\mathcal K})}+
\sum_{|\alpha|+j\le 9}\|\partial_x^j\Omega^\alpha
g\|_{L^2({\mathbb R}^3\backslash {\mathcal K})} \le \varepsilon,
\end{equation}
then \eqref{0.2} has a unique solution $u\in
C^\infty([0,T_\varepsilon]\times {\mathbb R}^3\backslash {\mathcal
K} )$, with
\begin{equation}\label{0.8}
T_\varepsilon = \exp(c/\varepsilon).
\end{equation}
\end{theorem}

We shall actually establish existence of limited regularity almost
global solutions $u$ for data $(f,g)\in H^9\times H^8$ satisfying
the relevant compatibility conditions and smallness assumptions
\eqref{dataob}. The fact then that $u$ must be smooth if $f$ and
$g$ are smooth and satisfy the compatibility conditions of
infinite order follows from standard local existence theorems (see
\S 9, \cite{KSS}).

As in \cite{KSS}, to prove this theorem it is convenient to show
that one can solve an equivalent nonlinear equation which has zero
initial data to avoid having to deal with issues regarding
compatibility conditions for the data.  We can then set up an
iteration argument for this new equation that is similar to the
one used in the proof of Theorem 3.3.

To make the reduction we first note that by local existence theory
(see, e.g., \cite{KSS}) if the data satisfies \eqref{dataob} with
$\varepsilon$ small we can find a local solution $u$ to $\square
u= Q(u')$ in $0<t<1$ that satisfies
\begin{multline}\label{0}
\sup_{0\le t\le 1}\sum_{|\alpha|\le 10} \Bigl( \|Z^\alpha u'(t,\,
\cdot\, )\|_{L^2({\mathbb R}^3\backslash {\mathcal K})}
\\
+\|(1+|x|)^{-1/2}Z^\alpha u'\|_{L^2(\{(s,x)\in [0,t]\times
{\mathbb R}^3\backslash {\mathcal K}) \})}\Bigr) \le C\varepsilon,
\end{multline}
for some uniform constant $C$.

Using this local solution we can set up our iteration.  We first
fix a  bump function $\eta\in C^\infty({\Bbb R})$ satisfying
$\eta(t)=1$ if $t\le 1/2$ and $\eta(t)=0$ if $t>1$.  If we set
$$u_0=\eta u$$
then
$$\square u_0=\eta Q(u')+[\square,\eta]u.$$
So $u$ will solve $\square u=Q(u')$ for $0<t<T_\varepsilon$ if and
only if $w=u-u_0$ solves
\begin{equation}\label{4.10}
\begin{cases}
\square w=(1-\eta)Q((u_0+w)')-[\square, \eta](u_0+w)
\\
w|_{\partial \mathcal{K}}=0
\\
w(0,x)=\partial_t w(0,x)=0
\end{cases}
\end{equation}
for $0<t<T_\varepsilon$.

We shall solve this equation by iteration.  We set $w_0=0$ and
then define $w_k$, $k=1,2,3,\dots$ recursively by requiring that
\begin{equation}\label{4.11}
\begin{cases}
\square w_k = (1 - \eta)Q((u_0+w_{k-1})')
-[\square,\eta](u_0+w_{k})
\\
w_k|_{\partial \obstacle}=0
\\
w_k(0,x)=\partial_tw_k(0,x)=0.
\end{cases}
\end{equation}

To proceed, we let
\begin{multline*}M_k(T)=
\sup_{0\le t\le T}\sum_{|\alpha|\le 10} \Bigl( \|Z^\alpha
w_k'(t,\,
\cdot\, )\|_2 \\
+(\ln(2+t))^{-1/2}\|(1+|x|)^{-1/2}Z^\alpha w_k'\|_{L^2(\{(s,x): \,
0\le s\le t\})}\Bigr).
\end{multline*}
Then, if we use \eqref{3.4}, \eqref{KSob} and \eqref{0}, we
conclude that there is a uniform constant $C_1$ so that
\begin{equation*}
M_k(T_\varepsilon)\le C_1\varepsilon + C_1\ln(2+T_\varepsilon)
(\varepsilon+M_{k-1}(T_\varepsilon))^2
+C_1(\varepsilon+M_{k-1}(T_\varepsilon))^2,
\end{equation*}
for some uniform constant $C_1$, if $\varepsilon$ is small.  Since
$M_0\equiv 0$, an induction argument implies that, if the constant
$c$ occurring in the definition of $T_\varepsilon$ is small then
\begin{equation}\label{3.17}
M_k(T_\varepsilon)\le 2C_1, \quad k=1,2,\dots, \end{equation}
 for
small $\varepsilon>0$.

If we let
\begin{multline*} A_k(T)=\sup_{0\le t\le T}\sum_{|\alpha|\le
10}\Bigl(\|Z^\alpha (u'_k-u'_{k-1})(t,\cd)\|_{L^2(\rt\backslash
\obstacle)}
\\
+(\ln(2+t))^{-1/2}\|(1+|x|)^{-1/2}Z^\alpha
(u'_k-u'_{k-1})\|_{L^2(\{(s,x): 0\le s\le t, x\in \rt\backslash
\obstacle \})}\Bigr),
\end{multline*}
then the preceding argument can be modified to show that
\begin{equation}\label{ak}
A_k(T_\varepsilon)\le \tfrac12 A_{k-1}(T_\varepsilon), \quad
k=1,2,\dots.
\end{equation}
Estimates \eqref{3.17} and \eqref{ak} imply Theorem
\ref{theorem0.1}. \qed

\newsection{Pointwise estimates}

To prove existence theorems for quasilinear wave equations we need
some pointwise estimates for solutions of inhomogeneous wave
equations, as well as some weighted Sobolev inequalities. To
describe the bounds for the wave equation, let us start out by
considering pointwise estimates for solutions of the inhomogeneous
wave equation in Minkowski space,
\begin{equation}\label{4.1}
\begin{cases}
(\partial_t^2-\Delta)w_0(t,x)=G(t,x),\quad (t,x)\in \R_+\times\R^3\\
w_0(0,x)=\partial_t w_0(0,x)=0.
\end{cases}
\end{equation}
If
$$L=t\partial_t + \langle x,\nabla_x\rangle$$
is the scaling operator, then in \cite{KSS3} the following
estimate was proved
\begin{equation}\label{4.2}
(1+t)|w_0(t,x)|\le C\sum_{\substack{|\alpha|+\mu\le 3 \\ \mu\le
1}}\int_0^t\int_{\rt} |L^\nu Z^\alpha G(s,y)|\, \frac{dy ds}{|y|}.
\end{equation}

Using this estimate and arguments from \S 2, we can obtain related
estimates for solutions of the inhomogeneous wave equation,
\begin{equation}\label{4.3}
\begin{cases}
(\partial_t^2-\Delta)w(t,x)=F(t,x),\quad (t,x)\in\R_+\times\ext\\
w(t,x)=0,\quad x\in\partial\mathcal{K}\\
w(t,x)=0, \quad t\le 0.
\end{cases}
\end{equation}
outside of obstacles satisfying Ikawa's local energy decay bounds
\eqref{1.6}. If we assume, as before, that
$\mathcal{K}\subset\{x\in\R^3\, :\, |x|<1\}$ and that
$\mathcal{K}$ satisfies \eqref{1.4} or \eqref{1.6}, then the
following pointwise estimate was proved in \cite{KSS3} and
\cite{MS}, respectively.
\begin{theorem}\label{theorem4.1}
Let $w$ be a solution to \eqref{4.3}, and suppose that the local
energy decay bounds \eqref{1.4} hold for $\mathcal{K}$.  Then,
\begin{multline}\label{4.4}
(1+t+|x|)|L^\nu Z^\alpha w(t,x)|\le C\int_0^t \int_\ext
\sum_{\substack{|\beta|+\mu \le |\alpha|+\nu+7\\\mu\le\nu +1}}
|L^\mu
Z^\beta F(s,y)|\:\frac{dy\:ds}{|y|}\\
+C\int_0^t \sum_{\substack{|\beta|+\mu\le |\alpha|+\nu+4\\
\mu\le\nu +1}}\|L^\mu \partial^\beta F(s,\cd)\|_{L^2(|y|<2)}\:ds.
\end{multline}
\end{theorem}

The estimate for non-trapping obstacles (in which case one can
take one less derivative in the right side of \eqref{4.4}) was
proved in \cite{KSS3}.  It was observed in \cite{MS} that the same
arguments will give \eqref{4.4} for obstacles satisfying Ikawa's
bounds \eqref{1.6}.

In \cite{MNS}, it was observed that one has the following
estimates for $w'$.

\begin{theorem}\label{theorem3.2}
Let  $w$ be a solution to \eqref{4.3}.  Suppose that $F(t,x)=0$
when $|x|>10 t$.  Then, if $|x|<t/10$ and $t>1$,
\begin{equation}\begin{split}\label{4.5}
(1+t+|x|)|L^\nu &Z^\alpha w'(t,x)|\le
C\sum_{\substack{\mu+|\beta|\le \nu+|\alpha|+3 \\ \mu\le \nu+1}
}\int_{t/100}^t \int_\ext
|L^{\mu}Z^{\beta}F'(s,y)|\:\frac{dy\:ds}{|y|}\\
&+ C \sup_{0\le s\le t}(1+s)\sum_{\substack{|\beta|+\mu\le
|\alpha|+4+\nu \\ \mu\le\nu}}\|L^\mu Z^\beta F(s,\cd)\|_\infty\\
&+ C\sup_{0\le s\le t} (1+s)\sum_{\substack{|\beta|+\mu\le
|\alpha|+\nu+6\\ \mu\le \nu}}\int_0^s
\int_{\substack{||y|-(s-\tau)|\le 10\\ |y|\le
(1000+\tau)/2}}|L^\mu Z^\beta F(\tau,y)|\:\frac{dy\:d\tau}{|y|}\\
&+C\sup_{0\le s\le t}\sum_{\substack{|\beta|+\mu\le |\alpha|+\nu+7\\
\mu\le\nu +1}}\int_{s/100}^s \int_{|y|\ge (1+\tau)/10}|L^\mu
Z^\beta F(\tau,y)|\:\frac{dy\:d\tau}{|y|}.
\end{split}\end{equation}
\end{theorem}

To prove either of these two estimates we realize that inequality
\eqref{4.2} yields
\begin{multline}\label{4.6}
(1+t)|L^\nu Z^\alpha w(t,x)|\le C \int_{0}^t
\int_{\ext}\sum_{\substack{|\beta|+\mu\le
|\alpha|+\nu+3\\\mu\le\nu +1}}|L^\mu Z^\beta
F(s,y)|\:\frac{dy\:ds}{|y|}
\\+C\sup_{|y|\le 2,0\le s\le t}(1+s)
\sum_{\substack{|\beta|+\mu\le \nu+|\alpha|+2\\ \mu\le \nu}}\|L^\mu
\partial^\beta w(s,\cd)\|_{L^2(|x|<2)}.
\end{multline}
The proof of \eqref{4.6} is exactly like that of Lemma 4.2 in
\cite{KSS3}. The last term in \eqref{4.6} can be estimated using
the local exponential decay of energy and the free space
estimates. This is the term that is responsible for the last term
in \eqref{4.4} and the last three terms in \eqref{4.5}.

As we mentioned before, we also need some weighted Sobolev
estimates.  The first is an exterior domain analog of results of
Klainerman-Sideris \cite{KS}.

\begin{lemma}\label{lemma4.3}  Suppose that $u(t,x)\in
C^\infty_0(\mathbb{R}\times \mathbb{R}^3 \backslash \mathcal{K})$
vanishes for $x\in \partial\mathcal{K}$.  Then if $|\alpha|=M$ and
$\nu$ are fixed
\begin{multline}\label{4.7}
\|\langle t-r\rangle L^\nu Z^\alpha \partial^2 u(t,\cd)\|_2\le
C\sum_{\substack{|\beta|+\mu \le M+\nu + 1 \\ \mu \le
\nu+1}}\|L^\mu Z^\beta u'(t,\cd)\|_2 \\
+C\sum_{\substack{|\beta|+\mu \le M+\nu \\ \mu\le \nu}} \|\langle
t+r\rangle L^\mu Z^\beta (\partial_t^2-\Delta) u(t,\cd)\|_2 +
C(1+t)\sum_{\mu\le \nu} \|L^\mu u'(t,\cd)\|_{L^2(|x|<2)}.
\end{multline}
\end{lemma}

The other such estimate that we need is an exterior domain analog
of an estimate of Hidano and Yokoyama \cite{HY}.
\begin{lemma}\label{lemma4.4}
Suppose that $u(t,x)\in C^\infty_0(\R\times\ext)$ vanishes for
$x\in\partial\mathcal{K}$.  Then
\begin{multline}\label{hy}
r^{1/2}\langle t-r\rangle |\partial L^\nu Z^\alpha  u(t,x)|\le
C\sum_{\substack{|\beta|+\mu\le |\alpha|+\nu+2\\ \mu\le
\nu+1}}\|L^\mu Z^\beta u'(t,\cd)\|_2 \\+
C\sum_{\substack{|\beta|+\mu\le |\alpha|+\nu+1\\\mu\le\nu}}
\|\langle t+r\rangle L^\mu Z^\beta (\partial_t^2-\Delta)
u(t,\cd)\|_2 +C(1+t)\sum_{\mu\le \nu} \|L^\mu
u'(t,\cd)\|_{L^\infty(|x|<2)}.
\end{multline}
\end{lemma}

\newsection{$L^2$ Estimates.}

In addition to the pointwise estimates, to prove global and almost
global existence results for quasilinear wave equations outside of
obstacles, we require certain energy-type estimates.  Since the
operators $\{Z\}$ and $L$ do not preserve the Dirichlet boundary
conditions, these are considerably more technical than the
estimates that are used for the Minkowski space setting, which
just follow from standard energy estimates and the fact that the
$Z$ operators commute with the D'Alembertian, while $[\square,
L]=2\square$.

The existence theorems involve possibly non-diagonal systems.
Because of this we are led to proving $L^2$ estimates for
solutions $u\in C^\infty (\Rplus\times\ext)$ of the Dirichlet-wave
equation
\begin{equation}\label{5.1}
\begin{cases}
\Box_\gamma u=F\\
u|_{\partial\mathcal{K}}=0\\
u|_{t=0}=f,\quad \partial_t u|_{t=0}=g
\end{cases}
\end{equation}
where
$$(\Box_\gamma u)^I=(\partial_t^2-c_I^2\Delta)u^I
+\sum_{J=1}^D\sum_{j,k=0}^3
\gamma^{IJ,jk}(t,x)\partial_j\partial_k u^J,\quad 1\le I\le D.$$
We shall assume that the $\gamma^{IJ,jk}$ satisfy the symmetry
conditions
\begin{equation}\label{5.2}
\gamma^{IJ,jk}=\gamma^{JI,jk}=\gamma^{IJ,kj}
\end{equation}
as well as the size condition
\begin{equation}\label{5.3}
\sum_{I,J=1}^D \sum_{j,k=0}^3 \|\gamma^{IJ,jk}(t,x)\|_{\infty}\le
\delta/(1+t),
\end{equation}
for $\delta$ sufficiently small (depending on the wave speeds).
 The energy estimate will involve bounds for the gradient of the
perturbation terms
$$\|\gamma'(t,\cd)\|_\infty = \sum_{I,J=1}^D\sum_{j,k,l=0}^3
\|\partial_l \gamma^{IJ,jk}(t,\cd)\|_\infty,$$ and the energy form
associated with $\Box_\gamma$, $e_0(u)=\sum_{I=1}^D e_0^I(u)$,
where
\begin{multline}\label{5.4}
e_0^I(u)=(\partial_0 u^I)^2+\sum_{k=1}^3 c_I^2(\partial_k u^I)^2
\\+2\sum_{J=1}^D\sum_{k=0}^3
\gamma^{IJ,0k}\partial_0u^I\partial_ku^J -
\sum_{J=1}^D\sum_{j,k=0}^3
\gamma^{IJ,jk}\partial_ju^I\partial_ku^J.
\end{multline}

The most basic estimate will lead to a bound for
$$E_M(t)=E_M(u)(t)=\int \sum_{j=0}^M e_0(\partial^j_t
u)(t,x)\:dx.$$

\begin{lemma}\label{lemma5.1}
Fix $M=0,1,2,\dots$, and assume that the perturbation terms
$\gamma^{IJ,jk}$ are as above.  Suppose also that $u\in C^\infty$
solves \eqref{5.1} and for every $t$, $u(t,x)=0$ for large $x$.
Then there is an absolute constant $C$ so that
\begin{equation}\label{5.5}
\partial_t E^{1/2}_M(t)\le C\sum_{j=0}^M \|\Box_\gamma
\partial_t^ju(t,\cd)\|_2+C\|\gamma'(t,\cd)\|_\infty E^{1/2}_M(t).
\end{equation}
\end{lemma}

This estimate is standard, and for this estimate one can weaken
\eqref{5.3} by replacing the right side with $\delta$ for
$\delta>0$ sufficiently small. It is important to note that there
is no ``loss" of derivatives here in \eqref{5.5}. On the other
hand, if we wish to prove bounds involving the $\{Z, L\}$
operators our techniques lead to estimates where there is an
additional local term which unfortunately involves a loss of one
derivative.  To be more specific, if we let
\begin{equation}\label{5.6}
Y_{N_0,\nu_0}(t)=\int \sum_{\substack{|\alpha|+\mu\le
N_0+\nu_0\\\mu\le\nu_0}}e_0(L^\mu Z^\alpha u)(t,x)\:dx.
\end{equation}
then, if \eqref{5.3} holds, we have
\begin{multline}\label{5.7}
\partial_t Y_{N_0,\nu_0}\le C Y^{1/2}_{N_0,\nu_0} \sum_{\substack{
|\alpha|+\mu\le N_0+\nu_0\\ \mu\le\nu_0}} \|\Box_\gamma L^\mu
Z^\alpha
u(t,\cd)\|_2 + C \|\gamma'(t,\cd)\|_\infty Y_{N_0,\nu_0} \\
+C \sum_{\substack{|\alpha|+\mu\le N_0+\nu_0+1\\ \mu\le \nu_0}}
\|L^\mu \partial^\alpha u'(s,\cd)\|^2_{L^2(|x|<1)}.
\end{multline}

In the arguments that are used to prove the existence theorems we
are able to handle the contributions of the last term in
\eqref{5.7} by using the following result from \cite{MS}.

\begin{lemma}\label{lemma5.2}
Suppose that \eqref{1.6} holds, and suppose that $u\in C^\infty$
solves \eqref{5.1} and satisfies $u(t,x)=0$ for $t<0$.  Then, for
fixed $N_0$ and $\nu_0$ and $t>2$,
\begin{multline}\label{5.8}
\sum_{\substack{|\alpha|+\mu\le N_0+\nu_0\\\mu\le\nu_0}} \int_0^t
\|L^\mu \partial^\alpha u'(s,\cd)\|_{L^2(|x|<2)}\:ds \\
\le C \sum_{\substack{|\alpha|+\mu\le N_0+\nu_0+1\\\mu\le\nu_0}}
\int_0^t \left(\int_0^s \|L^\mu \partial^\alpha \Box
u(\tau,\cd)\|_{L^2(||x|-(s-\tau)|<10)}\:d\tau\right)\:ds.
\end{multline}
\end{lemma}

These are the main $L^2$ estimates that are needed in the proof of
the existence results.  Using them and variations of the weighted
space-space time norms described in \S 3 that involve $L$ as well
as the operators $\{Z\}$ we can prove existence theorems for
certain quadratic, quasilinear systems of the form
\begin{equation}\label{5.9}
\begin{cases}
\Box u = Q(du, d^2u),\quad (t,x)\in \Rplus\times\ext\\
u(t,\cd)|_{\partial\mathcal{K}}=0\\
u(0,\cd)=f,\quad \partial_t u(0,\cd)=g.
\end{cases}
\end{equation}
Here
\begin{equation*}
\Box=(\Box_{c_1},\Box_{c_2}, \dots, \Box_{c_D})
\end{equation*}
is a vector-valued multiple speed D'Alembertian with
\begin{equation*}
\Box_{c_I}=\partial_t^2-c_I^2\Delta.
\end{equation*}
We will assume that the wave speeds $c_I$ are positive and
distinct.  This situation is referred to as the nonrelativistic
case.  Straightforward modifications of the argument give the more
general case where the various components are allowed to have the
same speed.  Also, $\Delta=\partial_1^2+\partial_2^2+\partial_3^2$
is the standard Laplacian.  Additionally, when convenient, we will
allow $x_0=t$ and $\partial_0=\partial_t$.

We shall assume that $Q(du,d^2u)$ is of the form
\begin{equation}\label{5.11}
Q^I(du,d^2u)=B^I(du)+\sum_{\substack{0\le j,k,l\le 3\\ 1\le J,K\le
D}}B^{IJ,jk}_{K,l}\partial_l u^K \partial_j\partial_k u^J, \quad
1\le I\le D
\end{equation}
where $B^I(du)$ is a quadratic form in the gradient of $u$ and
$B^{IJ,jk}_{K,l}$ are real constants satisfying the symmetry
conditions
\begin{equation}\label{5.12}
B^{IJ,jk}_{K,l}=B^{JI,jk}_{K,l}=B^{IJ,kj}_{K,l}.
\end{equation}
To obtain global existence, we shall also require that the
equations satisfy the following null condition which only involves
the self-interactions of each wave family.  That is, we require
that
\begin{equation}\label{5.13}
\sum_{0\le j,k,l\le 3}B^{JJ,jk}_{J,l}\xi_j\xi_k\xi_l=0
\quad\text{whenever}\quad
\frac{\xi_0^2}{c^2_J}-\xi_1^2-\xi_2^2-\xi_3^2=0,\quad J=1,\dots,D.
\end{equation}

To describe the null condition for the lower order terms, we
expand
$$B^I(du)=\sum_{\substack{1\le J,K\le D\\ 0\le j,k\le
3}}A^{I,jk}_{JK}\partial_j u^J \partial_k u^K.$$ We then require
that each component satisfy the similar null condition
\begin{equation}\label{5.13a}
\sum_{0\le j,k\le 3} A^{J,jk}_{JJ}\xi_j\xi_k=0 \quad\text{whenever
} \quad \frac{\xi_0^2}{c^2_J}-\xi_1^2-\xi_2^2-\xi_3^2=0, \quad
J=1,\dots, D.
\end{equation}
Thus, the null condition \eqref{5.13}-\eqref{5.13a} is one that
only involves interactions of components with the same wave speed.


We can now state the main result in \cite{MNS}:

\begin{theorem}\label{theorem1.1}
Let $\mathcal{K}$ be a fixed compact obstacle with smooth boundary
that satisfies \eqref{1.6}.  Assume that $Q(du,d^2u)$ and $\Box$
are as above and that $(f,g)\in C^\infty(\ext)$ satisfy the
compatibility conditions to infinite order.  Then there is a
constant $\varepsilon_0>0$, and an integer $N>0$ so that for all
$\varepsilon<\varepsilon_0$, if
\begin{equation}\label{1.11}
\sum_{|\alpha|\le N}\|<x>^{|\alpha|}\partial^\alpha_x
f\|_2+\sum_{|\alpha|\le N-1}\|<x>^{1+|\alpha|}\partial_x^\alpha
g\|_2 \le \varepsilon
\end{equation}
then \eqref{5.9} has a unique solution $u\in
C^\infty([0,\infty)\times\ext)$.
\end{theorem}

This result extended earlier ones of \cite{KSS} and \cite{MS}.  In
\cite{MS} a weaker theorem  was proved where instead of assuming
the null conditions \eqref{5.13} and \eqref{5.13a}, the authors
assumed that for every $I$ one has
\begin{equation*}
\sum_{0\le j,k,l\le 3}B^{IJ,jk}_{J,l}\xi_j\xi_k\xi_l=0
\quad\text{whenever}\quad
\frac{\xi_0^2}{c^2_J}-\xi_1^2-\xi_2^2-\xi_3^2=0,\quad J=1,\dots,D,
\end{equation*}
and
\begin{equation*}
\sum_{0\le j,k\le 3} A^{I,jk}_{JK}\xi_j\xi_k=0\quad \text{for all
}\quad \xi \in \R \times \R^3,  \quad 1\le J,K\le D.
\end{equation*}
respectively.

The nonrelativistic system satisfying the above null condition
that we study serves as a simplified model for the equations of
elasticity.  In Minkowski space, such equations were studied and
shown to have global solutions by Sideris-Tu \cite{Si3},
Agemi-Yokoyama \cite{AY}, and Kubota-Yokoyama \cite{KY}.

One can also, as in \cite{KSS3}, prove almost global existence for
solutions of equations of the form \eqref{5.9} that do not involve
null conditions.

\end{document}